\documentclass[12pt]{article}

\begin{document}

\title{Asymptotics of an empirical bridge of regression on induced order statistics}

\author{Artyom Kovalevskii\thanks{pandorra@ngs.ru, Novosibirsk State Technical University, 
Novosibirsk State University. The research was supported by RFBR grant 17-01-00683}}

\date{}

\maketitle


\begin{abstract}
We propose a class of tests for linear regression on concomitants (induced order statistics).
These tests are based on sequential sums of regression residuals. 
We self-center and self-normalize these sums. The resulting process is called an empirical bridge.
We prove weak convergence of the empirical bridge in uniform metrics to a centered Gaussian process.
The proposed tests are of chi-square type.

\end{abstract}


\section{Introduction}

There are few approaches to regression models testing.
An empirical fluctuation process function is implemented in R package by Zeileis at al. \cite{ref162}.
The process uses recursive regression residuals proposed by Brown, Durbin, Ewans \cite{ref0167}.

MacNeill \cite{ref13} studied linear regression for time series. He obtained limit
processes for sequences of partial sums of regression residuals.
Later Bischoff \cite{ref0165} showed that the MacNeill's theorem holds in
a more general setting. Aue
et al. \cite{ref1} introduced a new test for polynomial regression
functions which is analogous to the classical likelihood test. 
This approach is developed in \cite{ref101}, \cite{ref102}.

Stute
\cite{ref15} proposed a class of tests for one-parametric case.

Our approach is adopted specifically for regression on induced order statistics (concomitants).
These models arise in applications \cite{ref10}. Partial results are proposed in \cite{ref11}, \cite{ref12}.
The proofs use the theory of induced order statistics. 

David \cite{ref017} and Bhattacharya \cite{ref014} have introduced induced order statistics (concomitants)
simultaneously. The asymptotic theory was developed in 
in \cite{ref011}, \cite{ref012}, \cite{ref013}, \cite{ref015}, \cite{ref016}, \cite{ref018}, \cite{ref019}, 
\cite{ref0195}, \cite{ref02}, \cite{ref047}, 
\cite{ref05}, \cite{ref051}, \cite{ref053}, \cite{ref054}, \cite{ref056}, \cite{ref16}, \cite{ref161}.
Strong convergence to a corresponding Gaussian process can be proved by methods of 
\cite{ref052}, \cite{ref08}, \cite{ref131}, \cite{ref134}.

New contributions to the theory deal with extremal order statistics \cite{ref055},  \cite{ref133}, \cite{ref135}.

\section{Main result and corollaries}

Let $(\delta_i,\xi_i, \eta_i) =(\delta_i,\xi_{i1},\ldots, \xi_{im}, \eta_i)$ be independent and identically distributed 
random vector rows, $\delta_i$ has uniform distribution on $[0,1]$, $i=1,\ldots,n$. Random variables $\delta_1,\ldots,\delta_n$ 
are not observed. They will be used for ordering.

We assume a linear regression hypothesis $\eta_i=\xi_i \theta+e_i$. 
Here 
\[
\xi_i \theta\stackrel{def}{=}
\sum_{j=1}^m \xi_{ij}\theta_j,
\] 
$\{e_i\}_{i=1}^n$ and $\{(\delta_i, \xi_i)\}_{i=1}^n$ are independent, $\{e_i\}_{i=1}^n$ are i.i.d., ${\bf E}\, e_1 =0$, 
${\bf Var} \, e_1=\sigma^2>0$.

Vector $\theta=(\theta_1,\ldots,\theta_m)$ and constant $\sigma^2$ are unknown. 

We order rows by the first component, that is, we change rows $j$ and $k$ while $j<k$ and $\delta_{j}>\delta_{k}$.
The result is matrix $(U,X,Y)$ with rows $(U_i,{\bf X}_i, Y_i)=(U_i,X_{i1}, \ldots, X_{im}, Y_i)$. 
So $U_1<\ldots<U_n$ a.s. are order statistics from uniform distribution on $[0,1]$.
Elements of matrix $(X,Y)$ are concomitants. 
Let $\varepsilon_i=Y_i-{\bf X}_i \theta$. 
Note that $\varepsilon_i$ are i.i.d. and have the same distribution as $e_1$. 
Sequences $\{\varepsilon_i\}_{i=1}^n$ and $\{(U_i, X_i)\}_{i=1}^n$ are independent.
 
Let $\widehat{\theta}$ be LSE:
\[
\widehat{{\theta}}=({X}^T{X})^{-1}{X}^T{Y}.
\]

It does not depend on the order of rows.

Let
$h(x)={\bf E}\{\xi_{1} | \delta_1=x\}$ be conditional expectation,
$L(x)=\int\limits_{0}^{x}h (s)\,ds$  be {\it induced 
theoretical generalised Lorentz curve} (see [1]),
\[
b^2(x)={\bf E} \left((\xi_1-h(x))^T (\xi_1-h(x)) \ |\ \delta_1=x\right)
\]
be a matrix of conditional covariances.

Let $g_{ij}={\bf E}\xi_{1i}\xi_{1j}$, \ $G=(g_{ij})_{i,j=1}^m$,
Then $G=\int_0^1 (b^2(x)+h^T(x)h(x))\, dx$.

Let $\widehat{\varepsilon}_i=Y_i - X_i \widehat{\theta}$,
$
\widehat{\Delta}_k=\sum\limits_{i=1}^{k}\widehat{\varepsilon}_i$, $\widehat{\Delta}_0=0$.

Let $Z_n=\{Z_n(t), \, 0 \le t \le 1\}$ be a piecewise linear random function with nodes
\[
\left(\frac{k}{n}, \
\frac{\widehat{\Delta}_k}{\sigma\sqrt{n}}\right).
\]

We designate weak convergence in $C(0,1)$ with uniform metrics by $\Longrightarrow$.

{\bf Theorem 1} {\em If $G$ exists, ${\rm det} \, G \ne 0$, 
then $Z_n\Longrightarrow~Z$. Here $Z$ is a centered Gaussian process with covariation function }
\[
K(s,t)=\min(s,t)- L(s) G^{-1} L^T(t), \
~s, 
t\in[0,1].
\]

Let $Z_n^0$ be an empirical bridge (see \cite{ref10}, \cite{ref11}, \cite{ref12}):
\[
Z_n^0(t)=\frac{\sigma}{\widehat{\sigma}} (Z_n(t)-t Z_n(1)), \ \ 0 \le t \le 1,
\]
with $\widehat{\sigma}^2=\sum_{i=1}^n \widehat{\varepsilon}_i^2/n$.
Let $L^0(t)=L(t)-tL(1)$.

Let $L_{n,j}$ be an empirical induced generalised Lorentz curve:
\[
L_{n,j}(t)=\frac{1}{n} \sum_{i=1}^{[nt]} X_{ij},
\]
$L_n=(L_{n,1},\ldots, L_{n,m})$,  $L_n^0(t)=L_n(t)-t L_n(1)$.

{\bf Corollary 1} {\em Let assumptions of Theorem 1 be held. 

1) Then $Z_n^0\Longrightarrow~Z^0$, 
a centered Gaussian process with covariation function }
\[
K^0(s,t)=\min\{s,t\}- st - L^0(s) G^{-1} (L^0(t))^T, \
~s, 
t\in[0,1].
\]

{\em 2) Let $d\geq 1$ be integer, }
\[{\bf q}= (Z_n^0(1/(d+1)), \ldots, Z_n^0(d/(d+1)) ),
\]
$\widehat{g}_{ij}=\overline{X_i X_j}=\frac{1}{n} \sum_{k=1}^n X_{ki} X_{kj}$, $\widehat{G}=(\widehat{g}_{ij})_{i,j=1}^m$, 
\[\widehat{K}^0(s,t)=\min(s,t)- st - L_n^0(s))^T \widehat{G}^{-1} (L_n^0(t))^T,
\]
$Q=(\widehat{K}^0(i/(d+1), j/(d+1)))_{i,j=1}^d$. 
{\em Then 
${\bf q}Q^{-1}{\bf q}^T$ converges weakly to a chi-squared distribution with $d$ degrees of freedom.}

Note that ordering by $\xi_{i1}$, $i=1,\ldots,n$, can be viewed as ordering by $\delta_i$ with
$h_1(x)=F^{-1}_{\xi_{11}}(x)$ (quantile function). In this case $L_1(t)=\int_0^t F^{-1}_{\xi_{11}}(x) \, dx$.  

The next corollary is proved in \cite{ref11}.

{\bf Corollary 2} {\em Let $Y_i=\theta_1 X_{i1} + \varepsilon_i$, $i=1,\ldots,n$, $\theta_1 \in {\bf R}$, 
$(X_{11},\ldots,X_{n1})$ are order statistics of i.i.d. $(\xi_{11},\ldots,\xi_{n1})$, random variables
$(\varepsilon_1,\ldots,\varepsilon_n)$ are i.i.d. and independent of them, $0<{\bf E}\, \xi_{11}^2<\infty$,
${\bf E} \, \varepsilon_1=0$, $0<{\bf Var} \, \varepsilon_1=\sigma^2<\infty$. Then $Z_n \Rightarrow Z$, a
centered Gaussian process with covariance function }
\[\min(s,t)-L_1(s)L_1(t)/{\bf E}\, \xi_{11}^2.
\]

The next corollary is a partial case of Theorem 1 in \cite{ref12}.

{\bf Corollary 3} {\em Let $Y_i=\theta_1 X_{i1}+\theta_2 + \varepsilon_i$, $i=1,\ldots,n$, $\theta_1, \ \theta_2 \in {\bf R}$, 
$(X_{11},\ldots,X_{n1})$ are order statistics of i.i.d. $(\xi_{11},\ldots,\xi_{n1})$, random variables
$(\varepsilon_1,\ldots,\varepsilon_n)$ are i.i.d. and independent of them, $0<{\bf Var}\, \xi_{11}<\infty$,
${\bf E} \, \varepsilon_1=0$, $0<{\bf Var} \, \varepsilon_1=\sigma^2<\infty$. Then $Z_n \Rightarrow Z$, a
centered Gaussian process with covariance function}
\[
\min(s,t)-st - L_1^0(s) L_1^0(t)/{\bf Var}\, \xi_{11}.
\]

\section{Proof of Theorem 1}

Note that 
\[
\widehat{\Delta}_k=\sum_{i=1}^k (Y_i - X_i \widehat{\theta})
=
\sum_{i=1}^k (X_i(\theta- \widehat{\theta})+\varepsilon_i)
\]
\[
=
\sum_{i=1}^k (X_i(\theta- ({X}^T{X})^{-1}{X}^T{Y})+\varepsilon_i)
\]
\[
=
\sum_{i=1}^k (X_i(\theta- ({X}^T{X})^{-1}{X}^T(X\theta+\varepsilon))+\varepsilon_i)
\]
\[
=
\sum_{i=1}^k (\varepsilon_i-X_i({X}^T{X})^{-1}{X}^T\varepsilon).
\]

Note that $X_{[nt]}/n \to L(t)$ a.s. uniformely on compact sets,
and ${X}^T{X}/n \to G$ a.s.

So we study process
\[
\left\{
\sum_{i=1}^{[nt]} (\varepsilon_i-L(t) G^{-1}{X}^T\varepsilon), \ \
t\in [0,1]
\right\}.
\]

This process is a bounded linear functional of $(m+1)$-dimensional process
\[
\left\{
\sum_{i=1}^{[nt]} ({X}_i\varepsilon_i, \ \varepsilon_i), \ \
t\in [0,1]
\right\}.
\]

We use the functional central limit theorem for induced order statistics by Davydov and Egorov \cite{ref047}.

We assume that $\eta_i=\xi_i \theta+e_i$, $\{e_i\}_{i=1}^n$ and $\{\xi_i\}_{i=1}^n$ are independent, $\{e_i\}_{i=1}^n$ are i.i.d., ${\bf E}\, e_1 =0$, 
${\bf Var} \, e_1=\sigma^2>0$.

Let see rows $(\delta_i,\xi_i e_i, e_i) =(\delta_i,\xi_{i1}e_i,\ldots, \xi_{im}e_i, e_i)$. 
We have 
\[
{\bf E}(\xi_{1}e_1\ | \ \delta_1=x)=0, \ \ 
{\bf E}(e_1\ | \ \delta_1=x)=0, \ \ x \in [0,1].
\]

The conditional covariance matrix of the vector  $(\xi_1 e_1, e_1)$ is
\[
\tilde{b}^2(x)=
{\bf E} \left((\xi_1 e_1, e_1)^T (\xi_1 e_1, e_1)\ |\ \delta_1=x\right) = 
\sigma^2\left(
\begin{array}{cc} 
b^2(x)+h^T(x)h(x) & h^T(x) \\
h(x) & 1
\end{array}
\right).
\]

Let $\tilde{b}(x)$ be an upper triangular matrix such that $\tilde{b}(x)\tilde{b}(x)^T=\tilde{b}^2(x)$.
Then
\[
\tilde{b}(x)=
\sigma\left(
\begin{array}{cc} 
b(x) & h^T(x) \\
{\bf 0} & 1
\end{array}
\right).
\]

Here $b(x)$ is an upper triangular matrix such that ${b}(x){b}(x)^T={b}^2(x)$.
By Theorem 1 of Davydov and Egorov \cite{ref047} the process
\[
\left\{
\frac{1}{\sqrt{n}} 
\left(
\sum_{i=1}^{[nt]} X_i \varepsilon_i, \ \sum_{i=1}^{[nt]} \varepsilon_i
\right)^T, 
\ t \in[0,1]\right\}
\]
converges weakly in the uniform metrics to the Gaussian process
\[
\left\{
\int_0^t 
\tilde{b}(x) \, d{\bf W}_{m+1}(x), \ t \in[0,1]
\right\}.
\]

Here ${\bf W}_{m+1}=(W_1,\ldots, W_{m+1})^T$ is an $(m+1)$-dimensional standard Wiener process.

So process
\[
\left\{
\frac{1}{\sigma\sqrt{n}} 
\sum_{i=1}^{[nt]} (\varepsilon_i-L(t) G^{-1}{X}^T\varepsilon)^T, \ \
t\in [0,1]
\right\}
\]
converges weakly in the uniform metrics to the Gaussian process
$Z=\{Z(t), \ t \in [0,1]\}$,
\[
Z(t)=
W_{m+1} (t) - L(t) G^{-1} \left(
\int_0^1
b(x) \, d{\bf W}_{m}(x)+
\int_0^1 
h^T(x) \, d{W}_{m+1}(x) 
\right)^T,
\]

${\bf W}_{m}=(W_1,\ldots,W_m)^T$.

By the noted convergencies  $X_{[nt]}/n \to L(t)$ a.s. uniformely on compact sets,
${X}^T{X}/n \to G$ a.s., the process $Z_n$ has the same weak limit $Z$.

The covariance function of the limiting Gaussian process $Z$ is
\[
K(s,t)={\bf E} Z(s) Z(t) 
\]
\[
= \min(s,t) - L(s)G^{-1}\int_0^t h^T(x) \, dx -
L(t)G^{-1}\int_0^s h^T(x) \, dx 
\]
\[
+ L(s) G^{-1} \int_0^1 (b^2(x)+h^T(x)h(x))\, dx \, G^{-1} L^T(t)
\]
\[
=\min(s,t) - L(s) G^{-1} L^T(t).
\]

The proof is complete.

{\bf Acknowledgement}

The research was supported by RFBR grant 17-01-00683.

\end{document}